\documentclass[12pt]{article}
\pdfoutput=1
\usepackage[utf8]{inputenc}
\usepackage{fullpage}
\usepackage[small,bf]{caption}
\usepackage{float}
\usepackage{geometry}
\usepackage{verbatim}
\usepackage{listings}
\usepackage{graphicx}
\usepackage{algpseudocode}
\usepackage{enumerate, enumitem}
\usepackage{booktabs}
\usepackage{hyperref}
\hypersetup{colorlinks=true,
            urlcolor=blue,
            linkcolor=blue,
            citecolor=red}
\usepackage{url}
\usepackage{xcolor}
\usepackage{tikz}
\usetikzlibrary{math, positioning, calc, shapes, arrows.meta, patterns}
\usepackage{xparse}
\usepackage{amsmath,amsfonts}
\usepackage{subcaption}
\usepackage{textcomp}
\usepackage{caption}
\usepackage{threeparttable}
\usepackage{siunitx}

\definecolor{codegreen}{rgb}{0,0.6,0}
\definecolor{codegray}{rgb}{0.5,0.5,0.5}
\definecolor{codepurple}{rgb}{0.58,0,0.82}
\definecolor{backcolour}{rgb}{0.95,0.95,0.98}

\lstdefinestyle{mystyle}{ 
    backgroundcolor=\color{backcolour},   
    commentstyle=\color{codegreen}, 
    keywordstyle=\color{magenta},
    numberstyle=\tiny\color{codegray}, 
    stringstyle=\color{codepurple},
    basicstyle=\ttfamily\small, 
    breakatwhitespace=false,         
    breaklines=true,                 
    captionpos=t,                    
    keepspaces=true,                 
    numbers=left,                    
    numbersep=5pt,                  
    showspaces=false,                
    showstringspaces=false, 
    showtabs=false,                  
    tabsize=2}

\definecolor{seagreen}{rgb}{0.18, 0.55, 0.34}
\definecolor{mediumviolet-red}{rgb}{0.78, 0.08, 0.52}
\definecolor{khaki}{rgb}{0.94, 0.9, 0.55}

\lstdefinelanguage{mypython}
{
	keywords=[1]{from, import, as, assert, not, print, nonneg, boolean},
	keywordstyle=[1]{\color{mediumviolet-red}},
	keywords=[2]{cp, np, pd, cvxpy, numpy, pandas, sum, multiply, Variable, array,
    max, unique, sum_largest, hstack, matmul, Problem, Minimize, solve, range, set,
    max_daily_powers},
	keywordstyle=[2]{\color{seagreen}},
	upquote=true,
	showstringspaces=false,
	basicstyle=\ttfamily,
	columns=fullflexible,
	keepspaces=true,
	emph={True,False,def,return,float,class,match,switch,len},
	emphstyle={\color{seagreen}},
	belowskip=1em,
	aboveskip=1em,
    morecomment=[l]{\#}
}

\usepackage[
backend=bibtex,
style=numeric-comp,
bibstyle=ieee,
sorting=ynt,
firstinits=true,
url=false, doi=false, eprint=false,
]{biblatex}
\newcommand{\eg}{{\it e.g.}}

\newcommand{\BEQ}{\begin{equation}}
\newcommand{\EEQ}{\end{equation}}
\newcommand{\BEAS}{\begin{eqnarray*}}
\newcommand{\EEAS}{\end{eqnarray*}}

\newcommand{\reals}{{\mbox{\bf R}}}

\usepackage{relsize} 

\begin{document}

\title{Home Battery Dispatch under a Tiered Peak Power Tariff}

\author{David Pérez-Piñeiro\thanks{Department of Chemical Engineering,
Norwegian University of Science and Technology} \and Sigurd
Skogestad\footnotemark[1] \and Stephen Boyd\thanks{Department of Electrical
Engineering, Stanford University} }

\maketitle

\begin{abstract}
We consider the problem of operating a battery in a home connected to the
grid to minimize electricity cost, which combines an energy charge and a
tiered peak power charge based on the average of the $N$ largest daily
peak powers in each billing month. With perfect foresight of loads and prices,
the minimum cost is the solution of a mixed-integer linear program (MILP),
which provides a lower bound on the cost of any implementable policy. We
propose a model predictive control (MPC) policy that uses simple forecasts
of loads and prices and solves a small MILP at each time step. Numerical
experiments on one year of data from a home in Trondheim, Norway, show
that the MPC policy attains a cost within $1.7\%$ of the prescient bound,
and saves close to three times as much as the best rule-based policy we
consider.
\end{abstract}

\clearpage
\tableofcontents
\clearpage

\section{Introduction}
Retail electricity tariffs typically combine an energy charge, proportional
to consumption, with a peak power charge (also called a demand or capacity
charge) that reflects the customer's contribution to network
peaks~\cite{berg1983theory, henderson1983economics}. These charges have
traditionally applied to commercial and industrial customers, but as
residential electrification accelerates (heat pumps, electric vehicles,
on-site solar), several European countries have extended them to
households, with different designs now in use across the Nordic countries,
the Netherlands, and parts of central
Europe~\cite{wang2023dynamic, hofmann2025grid}.

Norway introduced such a tariff on 1~July 2022, motivated by the need to
align charges with network costs and defer grid
upgrades~\cite{verlo2020oppsummering}; capacity-based residential tariffs
in the Norwegian context have been studied in detail, see,
\eg,~\cite{bjarghov2022capacity}. Unlike a traditional demand charge
proportional to the single maximum power draw in a month, the Norwegian
tariff applies a \emph{tiered} charge to the \emph{average of the $N$
largest daily peak powers} (typically $N=3$). The monthly charge is
piecewise constant in this average, with jumps at a small number of
thresholds. In addition to this grid charge, the energy portion of the
bill combines a small time-of-use retail component with the Nord Pool
day-ahead wholesale price~\cite{NordPool2023}.

A home battery can shift consumption to low-price periods and shave peaks
to reduce the grid charge. These two objectives interact: greedy energy
arbitrage charges the battery at full rate whenever prices dip, driving
the grid draw into the highest peak tier. Minimizing total cost requires
jointly handling energy prices and the piecewise-constant peak cost.

\subsection{Our contribution}

We show how to formulate home battery dispatch under the Norwegian tiered
peak power tariff as a mixed-integer linear program (MILP). Solved with
perfect foresight of future loads and prices, it yields a \emph{prescient}
policy, whose cost lower bounds that of any implementable policy; solved on
a rolling horizon with forecasts, it yields an implementable \emph{model
predictive control} (MPC) policy. On one year of hourly data from a home in
Trondheim, the MPC policy comes within $1.7\%$ of the prescient lower bound,
and achieves close to three times the savings of the best rule-based policy
we consider.

\subsection{Related work}

\paragraph{Convex energy management.}
A unified framework for dynamic energy management based on convex
optimization appears in~\cite{moehle2019dynamic}, in which generators,
storage, and loads are modeled as a network of devices with convex
objectives and constraints. MPC variants that handle forecast uncertainty
through scenario-based robust optimization are developed
in~\cite{wytock2017dynamic}. Our formulation follows this modeling
approach; the only source of nonconvexity in our problem is the tiered
peak power charge, which introduces integer variables.

\paragraph{Battery scheduling under peak power charges.}
Most prior work on battery scheduling under peak power charges models the
charge as a linear penalty on the single largest power draw over the
billing period. The simplest methods are rule-based threshold policies,
which discharge the battery whenever the load exceeds a
target~\cite{karmiris2013peak}. More sophisticated methods use MPC to
jointly trade off peak reduction against energy cost, in both deterministic
and stochastic formulations~\cite{raoufat2018model, kumar2018stochastic,
risbeck2019economic}. Dynamic programming methods instead carry the running
peak in the state and use value-function approximation to manage the curse
of dimensionality~\cite{cholette2021battery}. The Norwegian tariff we
consider is different in two ways. The charge is based on the average of
the $N$ largest daily peaks rather than a single monthly maximum, and it is
tiered and piecewise constant rather than linear in that average. Because of
the averaging, a single running peak no longer summarizes the state, and the
tiered charge introduces the integer variables that make our problem
nonconvex.

\subsection{Outline}
In \S\ref{s-problem} we present a model of the home energy system. In
\S\ref{s-example} we introduce a running example based on data from a home
in Trondheim, Norway. In \S\ref{s-prescient} we formulate the prescient
problem as an MILP and establish a lower bound on achievable cost. In
\S\ref{s-MPC} we develop the MPC policy and compare it against the prescient
bound and three rule-based baselines. In \S\ref{s-forecasting} we describe
the methods used to generate the forecasts for MPC.
Appendix~\ref{s-sensitivity} presents a sensitivity analysis of MPC
performance to several problem parameters.

\section{Home energy system} \label{s-problem} 
\begin{figure}
\centering
\includegraphics[width=0.6\textwidth]{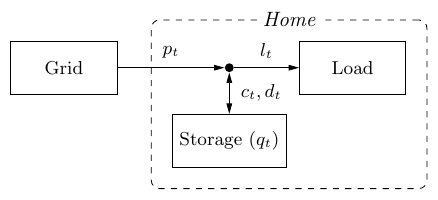}
\caption{Grid-connected home with storage. The power flows are the grid draw $p_t$, the charging and discharging powers $c_t$ and $d_t$, and the net load $l_t$. The storage device has an internal charge level $q_t$.}
\label{f-sys}
\end{figure}

We consider a grid-connected
home with a battery, as shown in figure~\ref{f-sys}.

\subsection{Power balance}
We use hourly time steps, indexed by $t=1,2, \ldots, T$. The load $l_t$ is
a net value that can include, for example, photovoltaic generation; it is
usually nonnegative. The storage power splits into a nonnegative charging
rate $c_t$ and a nonnegative discharging rate $d_t$. The grid power $p_t$
is the power drawn from the grid, with $0 \leq p_t \leq P$, where $P > 0$
is the maximum grid power. Power balance requires
\[
p_t + d_t - l_t - c_t = 0, \quad t=1, \ldots, T-1.
\]

\subsection{Storage device}
The storage device has a charge level $q_t$, with dynamics
\[
q_{t+1} = \eta_s q_{t} + h(\eta_c c_t - (1/\eta_d)d_t), \quad t=1, \ldots, T-1,
\]
where $\eta_s, \eta_c, \eta_d \in (0,1)$ are the storing, charging, and
discharging efficiencies, and $h$ is the time period in hours. The initial
charge level $q_1 = q_\mathrm{init}$ is given; we optionally impose a
terminal constraint $q_T = q_\mathrm{final}$. The charge level satisfies
\[
0 \leq q_t \leq Q, \quad t=1,\ldots, T,
\]
where $Q \geq 0$ is the storage capacity, and the charge and discharge
rates satisfy
\[
0 \leq c_t \leq C, \quad 0 \leq d_t \leq D, \quad t=1, \ldots, T-1,
\]
where $C$ and $D$ are the positive maximum rates.

\subsection{Cost function} \label{s-cost-function}

The cost function consists of an energy charge and a peak power charge.

\paragraph{Energy charge.}
The energy charge is $h\sum_{t=1}^{T-1} (\lambda_t^\mathrm{tou} +
\lambda_t^\mathrm{da})p_t$, where $\lambda_t^\mathrm{tou}$ are time-of-use
prices that vary by time of day and season, and $\lambda_t^\mathrm{da}$ are
day-ahead prices set daily through an auction market.

\paragraph{Peak power charge.}
The peak power charge is assessed monthly. We assume the time horizon spans $K$
whole months, denoted $k=1, \ldots, K$. Let $m_k$ denote the vector of daily
maximum powers in month $k$, and let $z_k$ denote the average of the $N$ largest
entries of $m_k$, where $N$ is a parameter (typically $N=3$). When $N=1$, $z_k$
is simply the maximum daily peak power over the month. We define the sum-largest
function $\psi(u,N)$ as the sum of the largest $N$ components of a vector $u$,
so $z_k = \psi(m_k,N)/N$.

The peak power cost in month $k$ is $\varphi(z_k)$, where $\varphi$ is
piecewise constant with the form
\BEQ \label{e-pk-pwr-cost}
\varphi(z) =
\begin{cases} \beta_1 & 0 \leq z \leq r_1\\
\beta_2 & r_1< z \leq r_2\\
\vdots \\
\beta_{L-1} & r_{L-2}< z \leq r_{L-1}\\
\beta_{L} & z > r_{L-1}
\end{cases}
\EEQ
where $0< \beta_1 < \beta_2 < \cdots
<\beta_{L}$ are the charges, and $0<r_1<r_2<\cdots <r_{L-1}$ are the thresholds.
When $\varphi(z_k) = \beta_j$, we say that the peak power cost is in tier~$j$.
We set $r_L=P$, so the last condition in \eqref{e-pk-pwr-cost} can be expressed
as $r_{L-1} < z \leq r_L$.

\paragraph{Total cost.}
The overall cost is
\[
h \sum_{t=1}^{T-1} (\lambda_t^\mathrm{tou} + \lambda_t^\mathrm{da}) p_t + \sum_{k=1}^K \varphi(z_k).
\]
The first term, the energy charge, is linear in the power values,
while the second term, the peak power charge, is nonconvex.

\subsection{Policies}
\paragraph{Information pattern.}
We assume that the following quantities are known.
\begin{itemize}
\item The system parameters $P$ (max grid power), $Q$ (battery capacity),
$C$ (max charge rate), $D$ (max discharge rate), and the efficiencies
$\eta_s$, $\eta_c$, and $\eta_d$.
\item The initial and final charge levels, $q_\mathrm{init}$ and
$q_\mathrm{final}$.
\item Time-of-use prices $\lambda_1^\mathrm{tou}, \ldots, \lambda_{T-1}^\mathrm{tou}$.
\item The parameter $N$ in the peak power calculation.
\item The costs $\beta_1, \ldots, \beta_{L}$ and thresholds $r_1, \ldots, r_{L-1}$
specifying the peak power cost tiers.
\end{itemize}

Day-ahead prices $\lambda_t^\mathrm{da}$ are announced daily at 13:00 for the
following day, so prices are known 12--35 hours in advance. The load $l_t$ is
known at the beginning of period $t$, but future loads are not known.

\paragraph{Prescient policy.}
All loads $l_1, \ldots, l_{T-1}$ and day-ahead prices
$\lambda_1^\mathrm{da}, \ldots, \lambda_{T-1}^\mathrm{da}$ are known in advance.
The charging and discharging powers are chosen to minimize total cost subject to
the constraints above. This policy is not implementable, but provides a lower
bound on the cost achievable by any implementable policy.

\paragraph{Implementable policy.}
An implementable policy respects the information pattern; in period $t$, it
chooses $c_t$ and $d_t$ based on loads and day-ahead prices
known at time $t$. Future values are not known, but can be forecast from past
observations.

\section{Running example} \label{s-example} We illustrate our methods on a
running example using real data from a home in Trondheim, Norway, over
2020--2022. We use 2020--2021 data to fit forecast models (see
\S\ref{s-forecasting}) and 2022 data to evaluate policies. The data and
source code are available at
\url{https://github.com/cvxgrp/home-battery-dispatch}.

\subsection{The Norwegian electricity tariff}
\label{s-norway}
A Norwegian residential electricity bill has three parts: a day-ahead
wholesale energy cost, a small retail energy markup, and a grid tariff that
includes the tiered peak power charge we study here.

\paragraph{Day-ahead prices.}
Wholesale electricity is traded on the Nord Pool spot market, which sets
hourly prices for each bidding zone one day in advance; Trondheim is in the
NO3 zone. The prices vary strongly over the day (morning and evening peaks,
a mid-day trough) and over the year, and rose sharply during the 2021--2022
European energy crisis, as shown in figure~\ref{f-da-prices}.

\paragraph{Retail energy charge.}
Retailers add a small time-of-use markup to the wholesale price, lower at
night and in winter. We use a markup representative of residential contracts
available in 2022.

\paragraph{Tiered peak power charge.}
Introduced in July 2022, Norway's residential capacity charge ties the
monthly charge to the customer's own peak demand~\cite{verlo2020oppsummering}.
Each distribution network operator sets a schedule of tiers, with a higher
charge when the monthly peak exceeds a threshold. We use the schedule of the
Trondheim network operator, shown in table~\ref{t-peak-tariff}. The charge
depends on $z_k$, the average of the $N = 3$ largest daily peaks in the
month, which reduces the effect of a single outlier day.

\begin{figure}
    \centering
    \includegraphics[width=.65\textwidth]{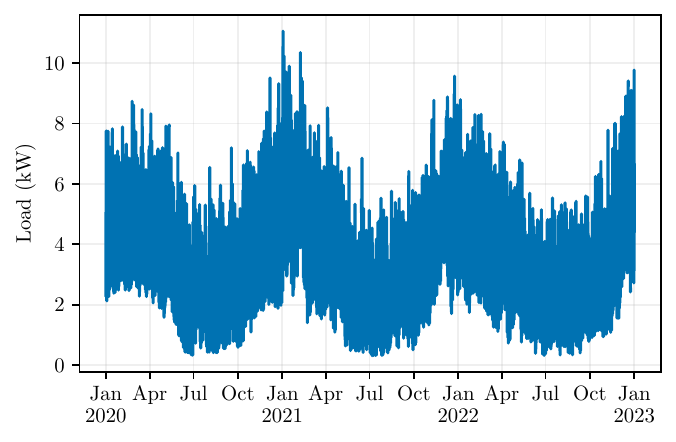}\\
    \includegraphics[width=.65\textwidth]{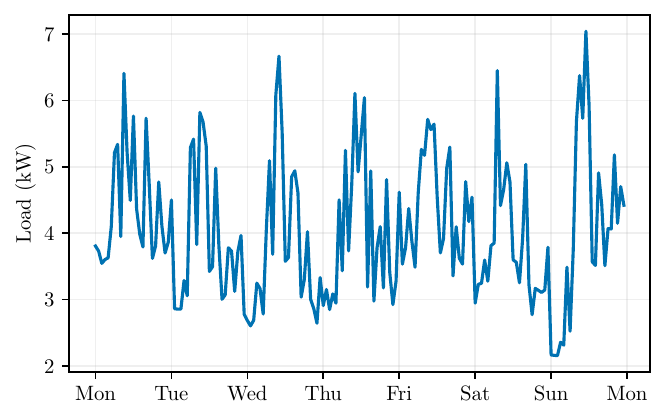}
    \caption{
    Hourly loads from a home in Trondheim, Norway.
    \emph{Top.}
    Three-year period 2020--2022.
    \emph{Bottom.}
    One week in January 2022.
    }
    \label{f-loads}
\end{figure}

\subsection{System parameters}
We use one-hour intervals, so $h=1$.
The hourly loads $l_t$ are shown in figure~\ref{f-loads}, over the full three
years (top) and one week in January 2022 (bottom).
The system parameters are
\[
P = 20\,\text{kW}, \quad Q=40\,\text{kWh}, \quad C=20\,\text{kW}, \quad D=20\,\text{kW},
\]
with efficiencies
\[
\eta_s = 0.99998, \quad \eta_c = 0.95, \quad \eta_d = 0.95.
\]
A complete charge or discharge takes 2~hours at full rate. The storing
efficiency $\eta_s$ corresponds to $1.5\%$ monthly self-discharge, typical
for lithium-ion batteries. The charging and discharging efficiencies each
represent a $5\%$ loss. We set the initial and final charge levels to
$q_\mathrm{init} = q_\mathrm{final} = Q/2$. In some
experiments we vary $Q$ around its nominal value of $40$\,kWh to study the
effect of storage capacity on cost savings.

\subsection{Price parameters}
\begin{figure}
    \centering
    \includegraphics[width=.65\textwidth]{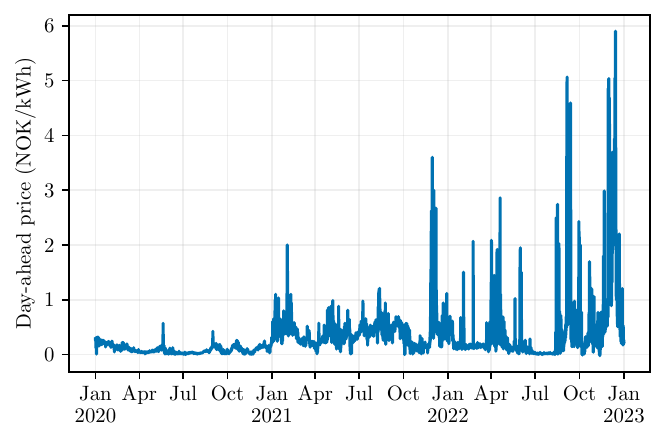}\\
    \includegraphics[width=.65\textwidth]{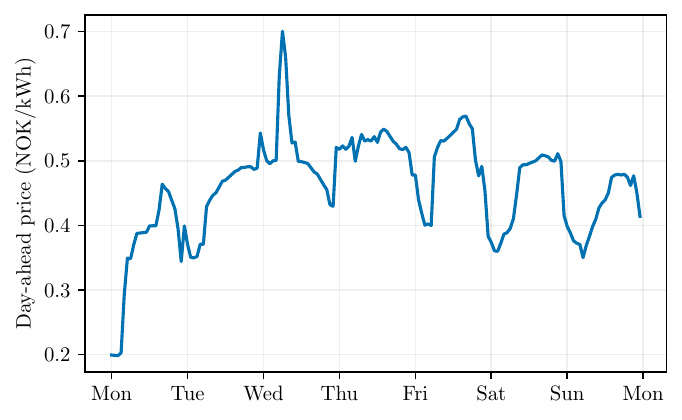}
    \caption{
    Day-ahead prices $\lambda_t^\mathrm{da}$ for Trondheim, Norway.
    \emph{Top.}
    Three-year period 2020--2022.
    \emph{Bottom.}
    One week in January 2022.
    }
    \label{f-da-prices}
\end{figure}

Figure~\ref{f-da-prices} shows day-ahead prices $\lambda_t^\mathrm{da}$
for Trondheim, in Norwegian Krone (NOK) per kWh, over the full three-year
period (top) and one week in January 2022 (bottom); data are sourced from
Nord Pool~\cite{NordPool2023}.

Table~\ref{t-tou-tariff} shows time-of-use prices, which vary by time of day and
season.

\begin{table}
    \centering
    \caption{Time-of-use prices $\lambda_t^\mathrm{tou}$ (NOK/kWh).}
    \label{t-tou-tariff}
    \begin{tabular}{lcc}
        \toprule
        & Day (hours 6--21) & Night (hours 22--5) \\
        \midrule
        Winter (Jan--Mar) & $0.30$ & $0.21$ \\
        Remaining months & $0.39$ & $0.30$ \\
        \bottomrule
    \end{tabular}
\end{table}

Table~\ref{t-peak-tariff} shows the tiered peak power cost. The monthly charge
is $\beta_l$ when $r_{l-1} < z_k \leq r_l$, where $z_k$ is the average of the
$N=3$ largest daily peak powers in month $k$.

\begin{table}
    \centering
    \caption{Tiered peak power cost (NOK/month).}
    \label{t-peak-tariff}
    \begin{tabular}{ccc}
        \toprule
        Tier $l$ & Threshold $r_l$ (kW) & Cost $\beta_l$ (NOK/month) \\
        \midrule
        1 & 2  & 83 \\
        2 & 5  & 147 \\
        3 & 10 & 252 \\
        4 & 15 & 371 \\
        5 & 20 & 490 \\
        \bottomrule
    \end{tabular}
\end{table}

\subsection{No storage baseline}
Without storage, the cost in 2022 is $25{,}052$ NOK, of which $22{,}028$ NOK
is the energy charge and $3{,}024$ NOK is the peak power charge.

\section{Prescient problem} \label{s-prescient}
In the prescient problem we assume the future loads and day-ahead prices
are known in advance. With $h=1$, the prescient problem is
\BEQ
\label{e-prescient}
\begin{array}{lll}
\mbox{minimize} & \sum_{t=1}^{T-1} (\lambda_t^\mathrm{tou} +
\lambda_t^\mathrm{da})p_t + \sum_{k=1}^K \varphi(z_k) \\
\mbox{subject to}
& p_t + d_t - l_t - c_t = 0, & t=1,\ldots,T-1 \\
& q_{t+1} = \eta_s q_t + \eta_c c_t - (1/\eta_d)d_t, & t=1,\ldots,T-1 \\
& q_1=Q/2, \quad q_T=Q/2 \\
& 0 \leq q_t \leq Q, & t=1,\ldots,T \\
& 0 \leq p_t \leq P, & t=1,\ldots,T-1 \\
& 0 \leq c_t \leq C, & t=1,\ldots,T-1 \\
& 0 \leq d_t \leq D, & t=1,\ldots,T-1
\end{array}
\EEQ
with variables $p, c, d \in \reals^{T-1}$ and $q \in \reals^{T}$, giving
about $4T$ scalar variables. Here $z_k = \psi(m_k,N)/N$, where $m_k$ is
the vector of daily maxima of $p_t$ in month $k$, and $\varphi$ is the
piecewise constant peak power cost function~\eqref{e-pk-pwr-cost}.

\subsection{MILP formulation}
The piecewise constant function $\varphi$ can be modeled by introducing
binary variables $s_{lk} \in \{0,1\}$ for tier $l=1,\ldots,L$ and
month $k=1,\ldots,K$, with
\[
\varphi(z_k) = \sum_{l=1}^L \beta_l s_{lk}, \qquad
z_k \leq \sum_{l=1}^L r_l s_{lk}, \qquad
\sum_{l=1}^L s_{lk} = 1.
\]
When $s_{lk}=1$, $z_k$ is assigned to tier $l$ with cost $\beta_l$.
This adds $LK$ binary variables.

The resulting problem is a mixed-integer convex problem (MICP); all
constraints are convex except the binary constraints $s_{lk} \in
\{0,1\}$. The objective is linear, most constraints are linear
equalities or inequalities, and the only nonlinearities involve $m_k$
(each entry is the maximum of a subset of the $p_t$, hence convex)
and $z_k = \psi(m_k, N)/N$, which is convex because $\psi(\cdot, N)$
is convex and nondecreasing~\cite[\S3.2.3]{boyd2004convex}.

The max and sum-largest functions are piecewise linear and can be
represented via linear inequalities, so~\eqref{e-prescient} can be
transformed to an equivalent MILP. Domain-specific languages (DSLs) for
convex optimization such as CVXPY~\cite{diamond2016cvxpy} automate these
transformations and generate an MILP that can be passed to solvers such
as Gurobi~\cite{gurobi} or open-source alternatives.

\subsection{CVXPY implementation}
We show how to formulate and solve problem~\eqref{e-prescient} using CVXPY.
We assume the parameters from \S\ref{s-example} and data vectors \verb|l|,
\verb|tou_prices|, \verb|da_prices|, and datetime index \verb|dt| are defined.

{\small
\begin{lstlisting}[language=mypython]
import cvxpy as cp
from pandas import unique

p = cp.Variable(T-1, nonneg=True)
c = cp.Variable(T-1, nonneg=True)
d = cp.Variable(T-1, nonneg=True)
q = cp.Variable(T, nonneg=True)
s = cp.Variable((K, L), boolean=True)

cons = [p + d - l - c == 0,
        q[1:] == eta_s*q[:-1] + eta_c*c - d/eta_d,
        q[0] == Q/2, q[-1] == Q/2,
        q <= Q, p <= P, c <= C, d <= D]

energy_cost = cp.sum(cp.multiply(tou_prices + da_prices, p))

peak_cost = 0
for k, month in enumerate(sorted(dt.month.unique())):
    days = unique(dt[dt.month == month].date)
    m_k = [cp.max(p[dt.date == day]) for day in days]
    z_k = cp.sum_largest(cp.hstack(m_k), N) / N
    peak_cost += beta @ s[k]
    cons += [z_k <= r @ s[k], cp.sum(s[k]) == 1]

prob = cp.Problem(cp.Minimize(energy_cost + peak_cost), cons)
prob.solve()
\end{lstlisting}
}
Lines 4--8 define variables matching the dimensions in~\eqref{e-prescient}.
Lines 10--13 encode constraints on power balance, storage dynamics, and
bounds. The loop (lines 18--23) builds the tiered peak power cost for
each month and adds the corresponding tier constraints, computing daily
maxima $m_k$ and the average of the $N$ largest via
\verb|cp.sum_largest|.

In roughly 20 lines we formulate and solve problem~\eqref{e-prescient}.
CVXPY automatically transforms expressions like \verb|cp.max| and
\verb|cp.sum_largest| into equivalent linear constraints before passing
the problem to an MILP solver.

\subsection{Running example}
We solve problem \eqref{e-prescient} using 2022 data with CVXPY and Gurobi.
After compilation, the problem has $35{,}783$ continuous variables and $60$
binary variables. On a laptop with Apple M1 Pro processor, it solves in about
$6$ seconds.

\paragraph{Costs.}
Table~\ref{t-costs-prescient} compares costs with and without storage. The
prescient cost is $21{,}204$ NOK, a savings of $15.4\%$ over the no-storage
baseline. Roughly two-thirds of the savings come from reduced energy charges.

\begin{table}
    \centering
    \caption{Annual electricity costs (NOK) for 2022 with 40\,kWh storage.}
    \label{t-costs-prescient}
    \begin{tabular}{lcccc}
        \toprule
        & Energy & Peak power & Total & Savings \\
        \midrule
        No storage & $22{,}028$ & $3{,}024$ & $25{,}052$ & --- \\
        Prescient  & $19{,}399$ & $1{,}805$ & $21{,}204$ & $15.4\%$ \\
        \bottomrule
    \end{tabular}
\end{table}

\begin{figure}
    \centering
    \includegraphics[width=.48\textwidth]{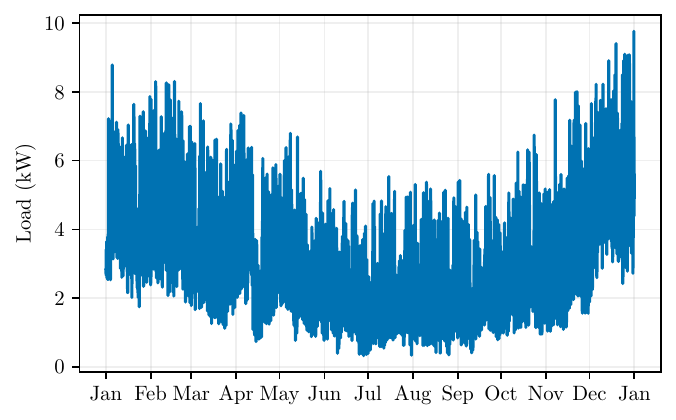}\hfill
    \includegraphics[width=.48\textwidth]{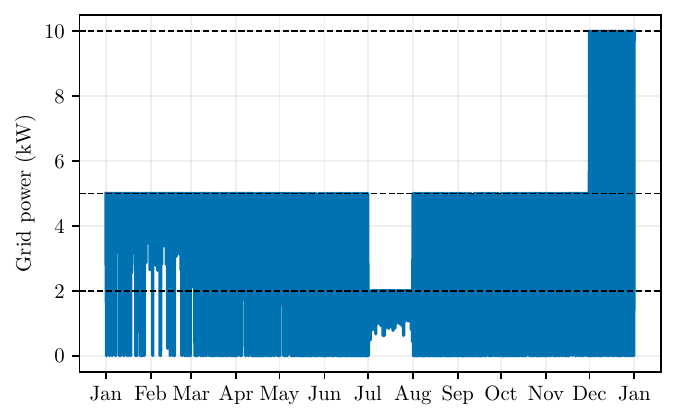}\\[1ex]
    \includegraphics[width=.48\textwidth]{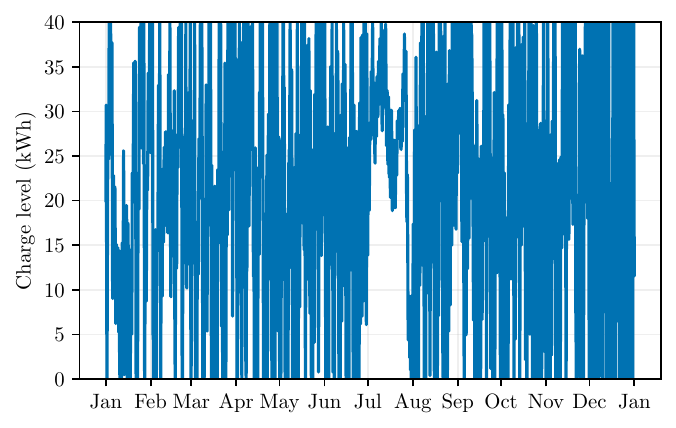}\hfill
    \includegraphics[width=.48\textwidth]{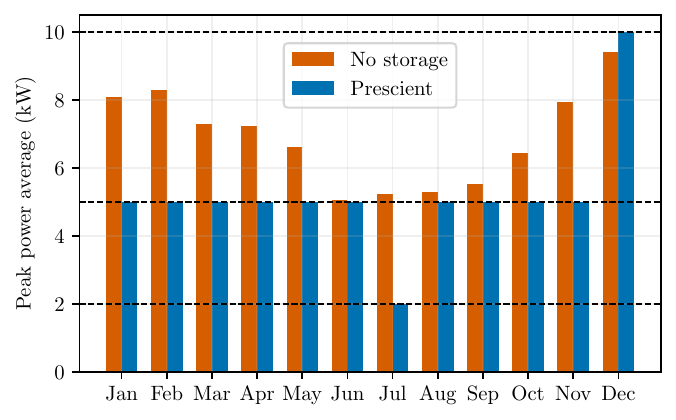}
    \caption{Prescient policy with 40\,kWh storage over 2022.
    Tier thresholds shown as dashed lines.
    \emph{Top left.} Load.
    \emph{Top right.} Grid power.
    \emph{Bottom left.} Charge level.
    \emph{Bottom right.} Average peak power $z_k$.}
    \label{f-flows-prescient}
\end{figure}

\begin{figure}
    \centering
    \includegraphics[width=.48\textwidth]{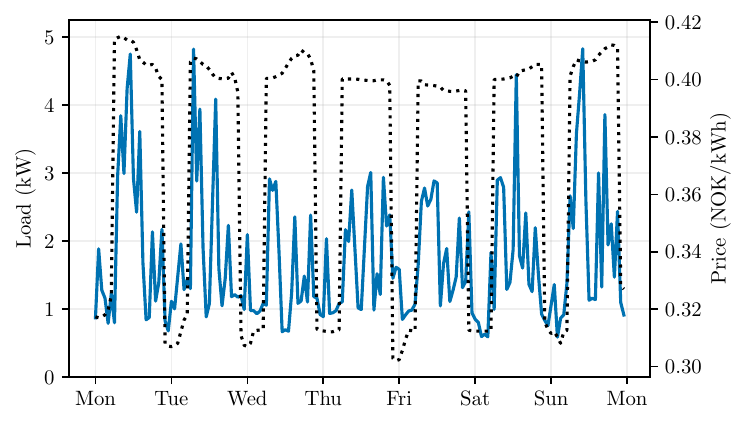}\hfill
    \includegraphics[width=.48\textwidth]{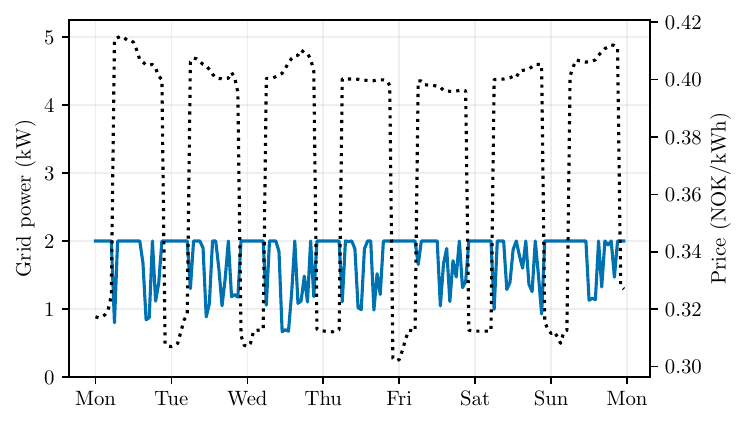}\\[1ex]
    \includegraphics[width=.48\textwidth]{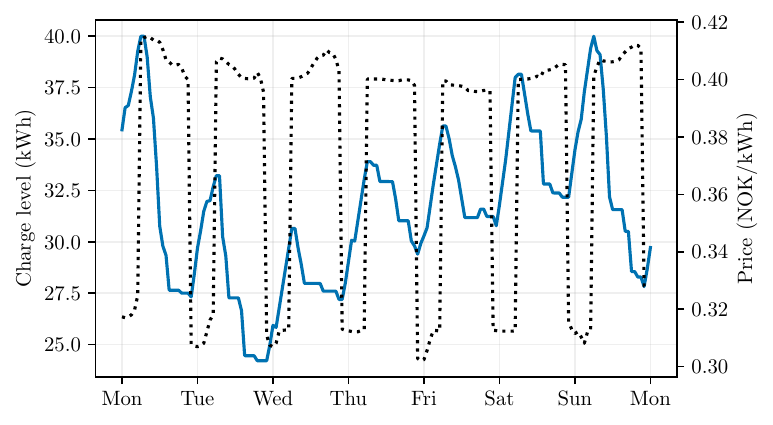}
    \caption{Prescient policy over one week in July 2022.
    Electricity prices $\lambda_t^\mathrm{tou} + \lambda_t^\mathrm{da}$
    shown as dotted lines.
    \emph{Top left.} Load.
    \emph{Top right.} Grid power.
    \emph{Bottom.} Charge level.}
    \label{f-week-prescient}
\end{figure}

\paragraph{Power flows.}
Figure~\ref{f-flows-prescient} shows the optimal power flows over the year. The
policy achieves tier~2 in ten months, tier~1 in July when loads are lowest, and
tier~3 only in December when winter peaks are highest. Figure~\ref{f-week-prescient}
shows one week in detail; the battery charges when prices are low and discharges
during load peaks, reducing both energy and peak power costs.

\begin{figure}
    \centering
    \includegraphics[width=.65\textwidth]{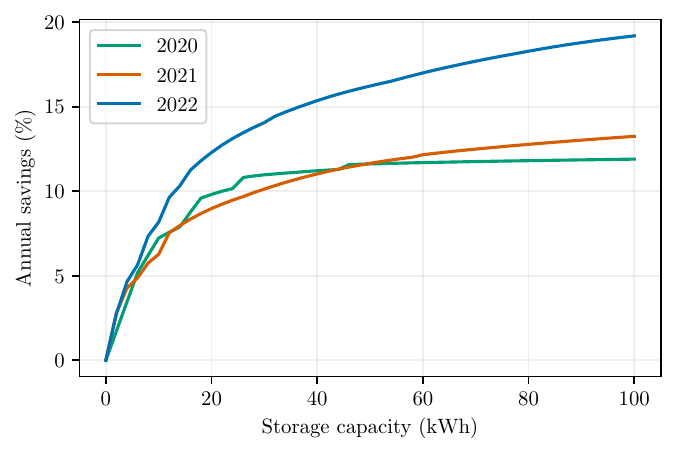}
    \caption{Annual savings versus storage capacity for 2020--2022.}
    \label{f-costs-vs-storage}
\end{figure}

\paragraph{Savings versus capacity.}
Figure~\ref{f-costs-vs-storage} shows annual savings versus storage capacity for 2020--2022. The curves exhibit diminishing returns; in 2022,
doubling capacity from $20$ to $40$ kWh increases savings from $12\%$ to
$15\%$. Without the binary tier constraints, the problem is an LP in $Q$, so the
curves are piecewise-linear and concave (a standard LP sensitivity
result); with them, tier transitions introduce small discontinuities but
leave the overall shape nearly concave.

\paragraph{Economic considerations.}
Figure~\ref{f-costs-vs-storage} lets us translate operational savings
into a break-even battery cost. Averaging across 2020--2022, a
$40$\,kWh battery saves about $69$\,NOK/kWh per year, so over a
ten-year service life the simple-payback break-even installed cost is
roughly $700$\,NOK/kWh. Residential lithium-ion systems in Norway
are installed at around $8{,}000$--$12{,}000$\,NOK/kWh before subsidy, so
a dedicated home battery does not yet pay back on operational savings
alone under this tariff. The policy is most valuable, then, for storage
that is already owned rather than bought for this purpose alone, such as
an electric-vehicle battery operated under
vehicle-to-grid~\cite{kempton2005vehicle}. For such assets the incremental
capital cost is negligible, and the savings we report are net.

\section{Model predictive control} \label{s-MPC}
The prescient problem provides a performance bound but requires perfect
foresight. We develop an implementable policy using model predictive control
(MPC). At each time step, we solve an optimization problem over a finite horizon
using forecasts, apply the first action, and repeat with updated information.

\subsection{MPC formulation}
At hour $t$, we observe the charge level $q_t$ and form forecasts of loads and
day-ahead prices over a planning horizon of $H$ hours. We let
$\hat{l}_{\tau|t}$ and $\hat{\lambda}_{\tau|t}^\mathrm{da}$ denote the forecasts
for period $\tau$ made at time $t$; for $\tau = t$ these are the known current
values. The forecasting method is described in \S\ref{s-forecasting}.
Let $k$ denote the current month and $M$ the number of months in the horizon.
The MPC problem is
\BEQ
\label{e-mpc}
\begin{array}{lll}
\mbox{minimize} & \sum_{\tau=t}^{t+H-1} (\lambda_\tau^\mathrm{tou} +
\hat{\lambda}_{\tau|t}^\mathrm{da}) p_\tau + \sum_{j=k}^{k+M-1} \varphi(z_{j|t}) \\
\mbox{subject to}
& p_\tau + d_\tau - \hat{l}_{\tau|t} - c_\tau = 0, & \tau = t,\ldots,t+H-1 \\
& q_{\tau+1} = \eta_s q_\tau + \eta_c c_\tau - (1/\eta_d)d_\tau, & \tau = t,\ldots,t+H-1 \\
& q_{t+H} = Q/2 \\
& 0 \leq q_\tau \leq Q, & \tau = t+1,\ldots,t+H \\
& 0 \leq p_\tau \leq P, & \tau = t,\ldots,t+H-1 \\
& 0 \leq c_\tau \leq C, & \tau = t,\ldots,t+H-1 \\
& 0 \leq d_\tau \leq D, & \tau = t,\ldots,t+H-1
\end{array}
\EEQ
with variables $p, c, d \in \reals^{H}$ and $q \in \reals^{H+1}$, where
the initial state $q_t$ is fixed at the observed value. The terminal
constraint $q_{t+H} = Q/2$ prevents myopic depletion of the battery.

The peak power terms require care because the monthly peak averages depend
on both past realized and future planned powers. Let $m_{j|t}$ be the
vector of daily maxima for month $j$ as evaluated at time $t$. For days
completed before $t$, the entry is the realized daily maximum (a
constant). For each day that intersects the horizon $[t, t+H-1]$, the
entry is the maximum of $p_\tau$ over the hours in that day, which may
combine past realized $p_\tau$ with decision variables $p_\tau$ for hours
$\tau \geq t$. The peak average $z_{j|t}$ is the average of the $N$
largest entries of $m_{j|t}$, or all entries if fewer than $N$ days are
represented.

\subsection{MILP formulation}
Problem \eqref{e-mpc} has the same structure as \eqref{e-prescient}, with
about $4H$ continuous variables. The piecewise constant $\varphi$ is
modeled using binary variables as in \S\ref{s-prescient}. A 30-day horizon
typically spans two months, giving $2L$ binary variables.

\subsection{LP enumeration}
An alternative to MILP is enumeration. Since each month selects exactly one
tier, we can solve \eqref{e-mpc} by enumerating over tier assignments. For each
assignment the tier variables are fixed and \eqref{e-mpc} reduces to a linear
program. With two months this gives $L^2$ linear programs; we solve them all and
select the minimum-cost solution. This is fast since $L$ is small, and is
useful when an MILP solver is unavailable.

\subsection{Running example}

\paragraph{Simple policies.}
We compare MPC against three rule-based policies that use only current
observations, without forecasts or optimization. All three are clipped to
respect charge-level and rate limits.

The \emph{peak shaving} policy discharges when the load exceeds a seasonal
target $\kappa_m$ and charges otherwise,
\[
(d_t, c_t) = \begin{cases}
(\min\{D,\, l_t - \kappa_m\},\; 0) & l_t > \kappa_m, \\
(0,\; \min\{C,\, \kappa_m - l_t\}) & l_t \leq \kappa_m,
\end{cases}
\]
where $\kappa_m$ depends on the current month. We use $\kappa_m = 10$\,kW
in winter (December--February), $\kappa_m = 5$\,kW in shoulder months
(March--May, September--November), and $\kappa_m = 2$\,kW in summer
(June--August), matching the upper boundary of tier~3, tier~2, and tier~1,
respectively. Loads are much higher in Norwegian winters than in summers.

The \emph{energy arbitrage} policy uses the current price to decide when to
charge and discharge. Let $\lambda_t = \lambda_t^\mathrm{tou} +
\lambda_t^\mathrm{da}$ be the retail price at hour $t$, and
$\bar{\lambda}_t$ its median over the calendar day that contains $t$.
The policy charges at full rate when the current price is below the median
and discharges against the load otherwise,
\[
(d_t, c_t) = \begin{cases}
(0,\; \min\{C,\, (Q - q_t)/\eta_c\})
 & \lambda_t < \bar{\lambda}_t, \\
(\min\{D,\, l_t\},\; 0)
 & \lambda_t \geq \bar{\lambda}_t.
\end{cases}
\]
The threshold is the daily median rather than a fixed hour schedule, so the
policy tracks the prevailing day-ahead price curve.

The \emph{capped arbitrage} policy merges the two by capping the grid draw
at the seasonal target $\kappa_m$ during charging,
\[
(d_t, c_t) = \begin{cases}
(0,\; \min\{C,\, (Q - q_t)/\eta_c,\, \kappa_m - l_t\})
 & \lambda_t < \bar{\lambda}_t, \\
(\min\{D,\, l_t\},\; 0)
 & \lambda_t \geq \bar{\lambda}_t.
\end{cases}
\]
Charging is thus clipped both to the charge level and to the tier cap,
so the policy exploits energy arbitrage when prices are low but never
drives the grid draw above its seasonal tier.

\paragraph{MPC parameters.}
We run MPC with horizon $H = 720$ hours (30 days) and $N = 3$, matching the
tariff definition, using the forecasts described in \S\ref{s-forecasting}. With
$L=5$ tiers, the horizon spans two months, giving $2L = 10$ binary variables.
Solving the MILP with Gurobi, each step takes around $0.2$ seconds.
Sensitivity to $H$, $N$, and the forecast method is reported in
Appendix~\ref{s-sensitivity}.

\paragraph{Costs.}
Table~\ref{t-mpc-results} compares all policies. MPC achieves a total cost of
$21{,}568$ NOK, a savings of $13.9\%$ relative to no storage and within $1.7\%$
of the prescient bound.
Peak shaving achieves $3.3\%$ savings by reducing peak power charges in
shoulder and summer months but forgoes energy arbitrage; even with a seasonal
target, the policy is sometimes too aggressive in summer (the battery
depletes serving routine load) and too passive in winter (loads rarely
exceed the target, so the battery stays idle). Energy arbitrage reduces
energy cost by nearly a thousand NOK but, by charging at full rate
whenever prices fall below the daily median, drives peak power into
tier~5 every month, giving a total cost \emph{higher} than no storage. Capped arbitrage, which clips the grid
draw during low-price charging to the seasonal target, avoids this failure
and achieves $5.0\%$ savings.

The gap between the best rule-based policy ($5.0\%$) and MPC ($13.9\%$) is
almost a factor of three, and quantifies the value of forecasting and joint
optimization. To isolate the contribution of modeling the peak power charge,
we also report an MPC variant that ignores it (energy-only). This policy
achieves the lowest energy cost but, like energy arbitrage, hits tier~5 every
month, yielding savings of only $1.6\%$. The jump from $1.6\%$ to $13.9\%$
shows that most of the value of MPC comes from jointly handling energy and
peak charges; arbitrage alone is not enough.

\begin{table}
    \centering
    \caption{Comparison of policies for 2022 with 40\,kWh storage. Costs in NOK.}
    \label{t-mpc-results}
    \begin{tabular}{lcccc}
        \toprule
        & Energy & Peak power & Total & Savings \\
        \midrule
        No storage        & $22{,}028$ & $3{,}024$ & $25{,}052$ & --- \\
        Peak shaving      & $22{,}009$ & $2{,}225$ & $24{,}234$ & $3.3\%$ \\
        Energy arbitrage  & $21{,}088$ & $5{,}880$ & $26{,}968$ & $-7.7\%$ \\
        Capped arbitrage  & $20{,}779$ & $3{,}024$ & $23{,}803$ & $5.0\%$ \\
        MPC (energy-only) & $18{,}774$ & $5{,}880$ & $24{,}654$ & $1.6\%$ \\
        MPC               & $19{,}279$ & $2{,}289$ & $21{,}568$ & $13.9\%$ \\
        Prescient         & $19{,}399$ & $1{,}805$ & $21{,}204$ & $15.4\%$ \\
        \bottomrule
    \end{tabular}
\end{table}

\begin{figure}
    \centering
    \includegraphics[width=.48\textwidth]{figures/load_year.pdf}\hfill
    \includegraphics[width=.48\textwidth]{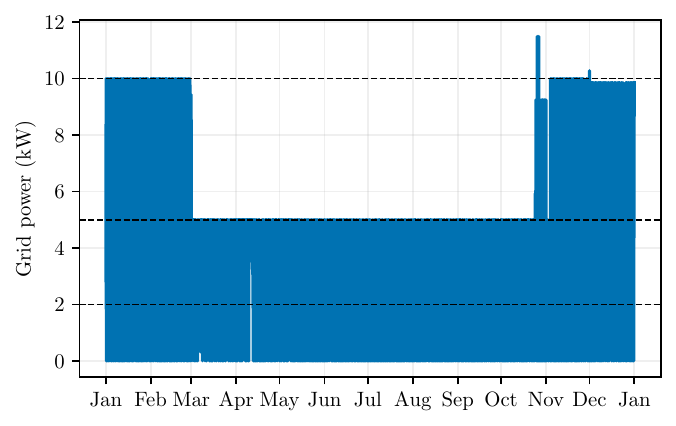}\\[1ex]
    \includegraphics[width=.48\textwidth]{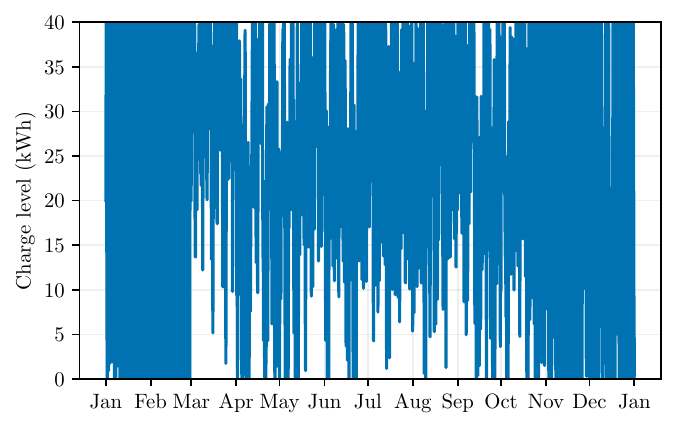}\hfill
    \includegraphics[width=.48\textwidth]{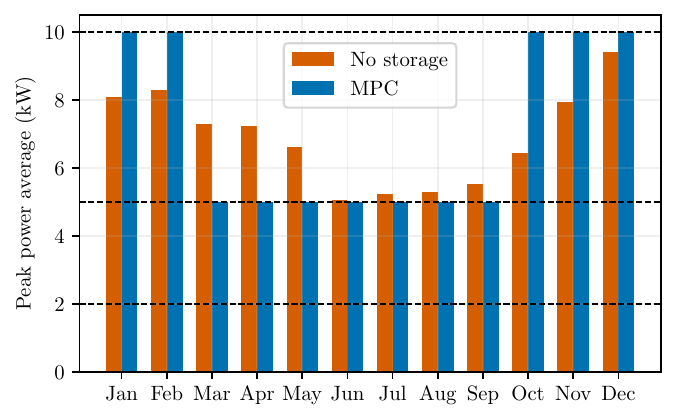}
    \caption{MPC policy with 40\,kWh storage over 2022.
    Tier thresholds shown as dashed lines.
    \emph{Top left.} Load.
    \emph{Top right.} Grid power.
    \emph{Bottom left.} Charge level.
    \emph{Bottom right.} Average peak power $z_k$.}
    \label{f-flows-mpc}
\end{figure}

\begin{figure}
    \centering
    \includegraphics[width=.48\textwidth]{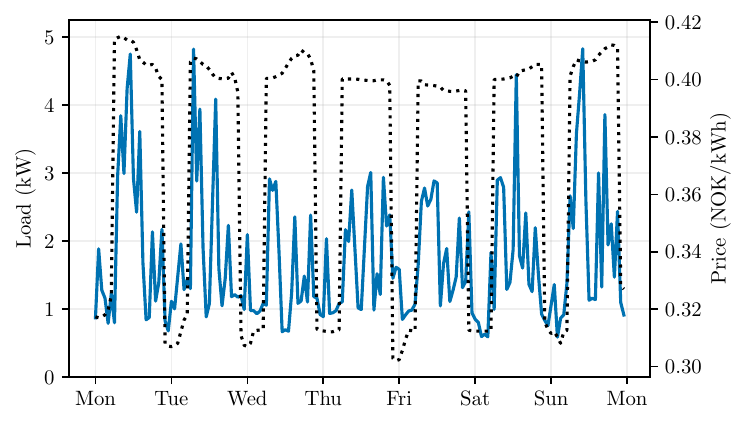}\hfill
    \includegraphics[width=.48\textwidth]{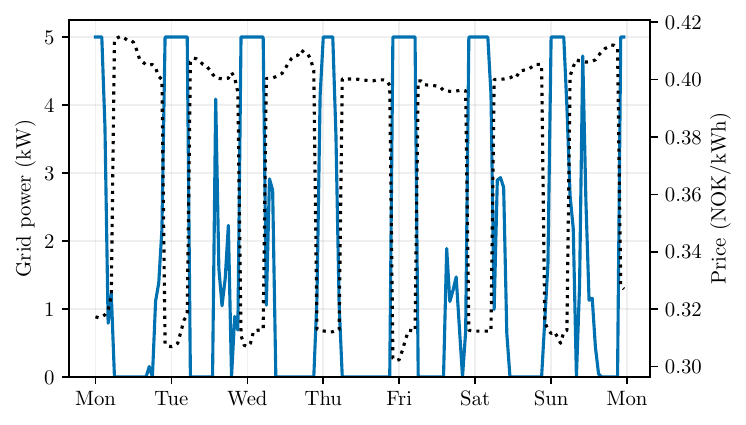}\\[1ex]
    \includegraphics[width=.48\textwidth]{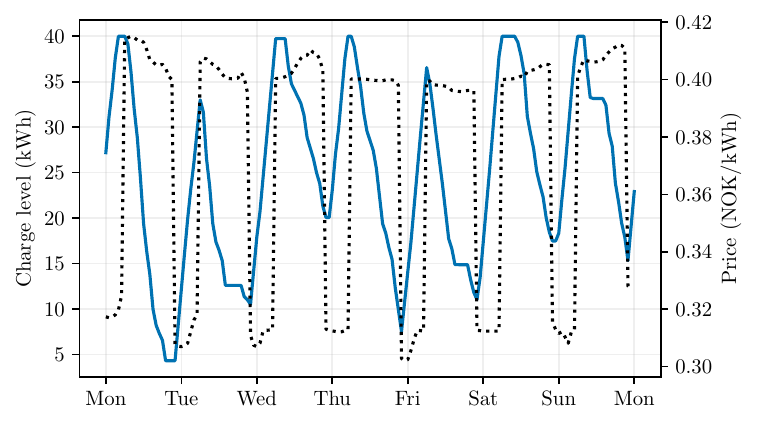}
    \caption{MPC policy over one week in July 2022.
    Electricity prices $\lambda_t^\mathrm{tou} + \lambda_t^\mathrm{da}$
    shown as dotted lines.
    \emph{Top left.} Load.
    \emph{Top right.} Grid power.
    \emph{Bottom.} Charge level.}
    \label{f-week-mpc}
\end{figure}

\paragraph{Power flows.}
Figure~\ref{f-flows-mpc} shows the MPC power flows over the year. MPC
achieves tier~2 in fewer months than prescient, since forecast errors
make it harder to anticipate the highest load days. The weekly detail
in figure~\ref{f-week-mpc} shows similar coordination to prescient, with
timing differences due to uncertainty in load and price forecasts.

\section{Forecasting} \label{s-forecasting}
In this section we describe a simple method to forecast a scalar time series
$x_1, x_2, \ldots$ using historical data, following
\textcite[Appendix~A]{moehle2019dynamic}.

\subsection{The baseline-residual forecast}
The baseline-residual forecast is given by
\[
\hat{x}_{\tau|t} = b_{\tau} + \hat{r}_{\tau|t},
\]
where $\hat{x}_{\tau|t} \in \reals$ is the prediction of quantity $x$ at
time $\tau \geq t$, made at time $t$.
The forecast is the sum of a seasonal baseline $b_\tau$, capturing periodic
patterns (diurnal, weekly, annual), and an autoregressive residual
$\hat{r}_{\tau|t}$ that accounts for short-term deviations from the baseline.
The baseline depends only on $\tau$; the residual depends on both $\tau$ and the
time $t$ at which the forecast is made.

\paragraph{Baseline component.} A simple model for the baseline is a sum of $K$
sinusoids,
\[
b_t = \beta_0 + \sum_{k=1}^K \left( \alpha_k \sin\left(\frac{2\pi t}{P_k}\right)
+ \beta_k \cos\left(\frac{2\pi t}{P_k}\right) \right),
\]
where $\alpha_k$ and $\beta_k$ are the coefficients and $P_k$ are the periods.  
In the usual case of Fourier series, the periods have the form $P_k=P/k$,
where $P$ is the fundamental period. Here we include terms for daily, weekly,
and annual variation.

To fit the $2K+1$ coefficients $\beta_0, \alpha_1, \beta_1, \ldots, \alpha_K,
\beta_K$ to historical training data $x_1, \ldots, x_T$, we use a pinball
(quantile) loss with $\ell_2$ regularization. The pinball loss for quantile $\eta
\in [0,1]$ is
\[
L_{\eta}(u) = \max\{\eta u, (\eta - 1) u\} = (\eta - 1/2) u + (1/2) |u|.
\]
Being piecewise linear, it is less sensitive to outliers than squared error.
Choosing $\eta < 0.5$ biases forecasts toward overestimation, while $\eta >
0.5$ biases toward underestimation. We find the coefficients by minimizing
\[
\sum_{t=1}^T L_{\eta}(b_t - x_t) + \lambda
\sum_{k=1}^K \nu_k (\alpha_k^2 + \beta_k^2),
\]
where $\lambda > 0$ is the regularization parameter and $\nu_k = k^2$ penalizes
higher harmonics more heavily. The $\ell_2$ penalty shrinks coefficients without
zeroing them out, since all harmonics contribute to a smooth baseline. The problem is convex and readily solved; good
values for $\eta$ and $\lambda$ can be chosen by cross-validation
\cite[\S7.10]{hastie2001elements}, \cite[\S13.2]{boyd2018introduction}.

\paragraph{Residual component.}
The residuals are $r_t = x_t - b_t$, $t = 1, \ldots, T$. We fit an
autoregressive (AR) model to predict the next $L$ residuals from the previous
$M$,
\[
(\hat{r}_{t+1|t}, \ldots, \hat{r}_{t+L|t}) = \Gamma(r_{t-M+1}, \ldots, r_t),
\]
where $\Gamma \in \reals^{L \times M}$ is the parameter matrix. We find
$\Gamma$ by minimizing
\[
\sum_{t=1}^T \sum_{\tau=t+1}^{t+L} L_{\eta}(\hat{r}_{\tau|t} - r_\tau) +
\lambda \|\Gamma\|_F^2,
\]
where $\|\cdot\|_F$ is the Frobenius norm; good values for $\eta$ and $\lambda$
can be chosen by cross-validation.

\subsection{Running example}

\begin{figure}
    \centering
    \includegraphics[width=.65\textwidth]{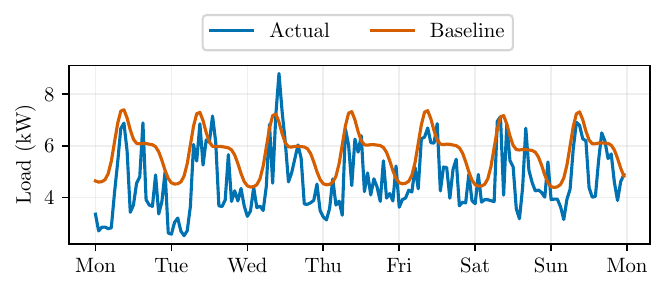}\\
    \includegraphics[width=.65\textwidth]{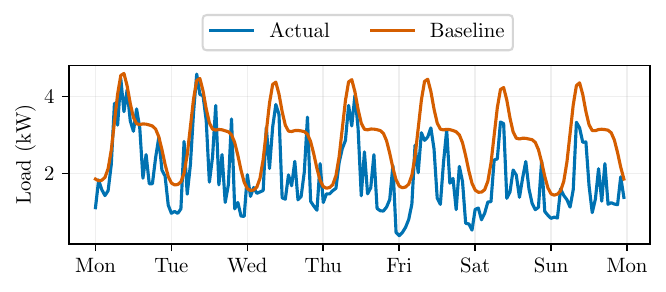}
    \caption{Baseline component (orange) and actual load (blue).
    \emph{Top.} One week in January 2022. \emph{Bottom.} One week in June 2022.}
    \label{f-load-baseline}
\end{figure}

\paragraph{Load forecasting.} We apply the baseline-residual method to load.
The baseline captures diurnal (24h), weekly, and annual periodicities with 4
harmonics each, giving periods
\[
\begin{array}{l}
P_1=24/1,~P_2=24/2,~P_3=24/3,~P_4=24/4,\\
P_5=168/1,~P_6=168/2,~P_7=168/3,~P_8=168/4,\\
P_9=8760/1,~P_{10}=8760/2,~P_{11}=8760/3,~P_{12}=8760/4.
\end{array}
\]
We fit these 25 parameters on two years of hourly data (2020--2021).
Figure~\ref{f-load-baseline} shows the baseline component against actual
load for one week each in January and June 2022.

\begin{figure}
    \centering
    \includegraphics[width=.65\textwidth]{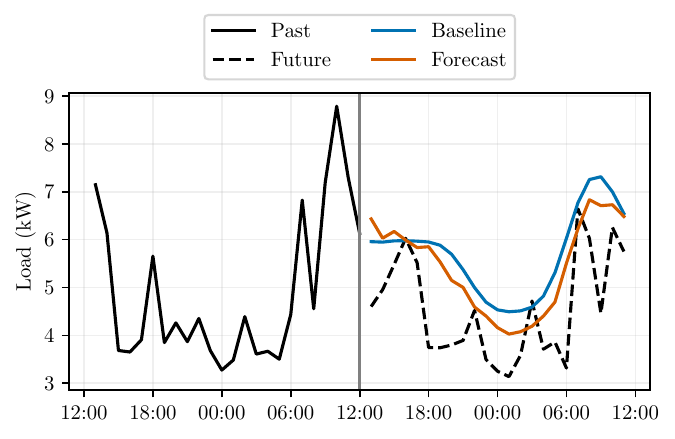}
    \caption{Load forecasts for a test day in May 2022, with vertical line marking
    the forecast hour. Black solid and dashed lines show realized and future
    load; blue shows the baseline component; orange the forecast (baseline + AR).}
    \label{f-load-forecast-comparison}
\end{figure}

For the residual component, we fit an AR model to predict residuals over the
next 23 hours from the previous 24. The $23 \times 24$ parameter matrix
$\Gamma$ is fit using quantile regression with $\ell_2$ regularization on the
same training period. Figure~\ref{f-load-forecast-comparison} compares the
baseline component alone to the full forecast (baseline + AR) for a test day in
May 2022. 

\begin{figure}
    \centering
    \includegraphics[width=.65\textwidth]{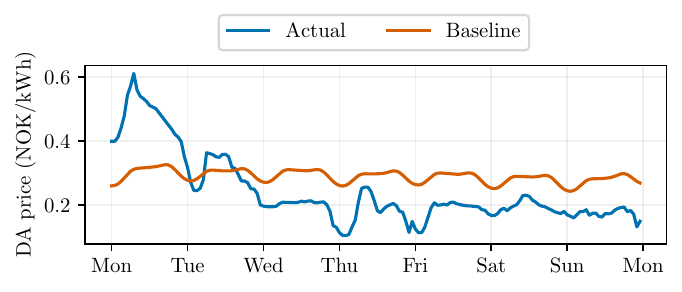}\\
    \includegraphics[width=.65\textwidth]{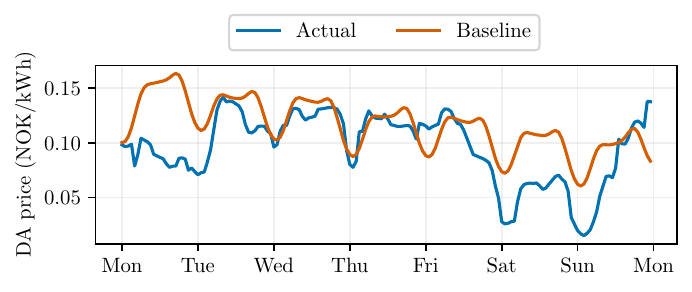}
    \caption{Baseline component (orange) and actual day-ahead prices (blue).
    \emph{Top.} One week in January 2022. \emph{Bottom.} One week in June 2022.}
    \label{f-price-baseline}
\end{figure}

\begin{figure}
    \centering
    \includegraphics[width=.65\textwidth]{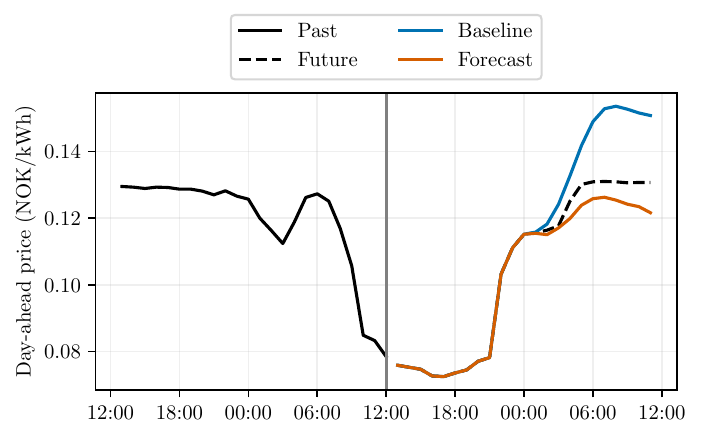}
    \caption{Price forecasts for a test day in May 2022, with vertical line marking
    the forecast hour. Black solid and dashed lines show realized and future
    prices; blue shows the baseline component; orange the forecast (baseline + AR).}
    \label{f-price-forecast-comparison}
\end{figure}

\paragraph{Forecasting day-ahead prices.} We apply the same method to day-ahead
prices. Since prices are announced at 13:00 for the next day, the known horizon
varies from 12 to 35 hours depending on the time of day.
Figure~\ref{f-price-baseline} shows the baseline component against
actual prices for one week each in January and June 2022.
Figure~\ref{f-price-forecast-comparison} compares the baseline component alone
to the full forecast (baseline + AR) for a test day in May 2022. 

\section{Conclusions}
\label{s-conclusions}

We have developed an MPC policy for operating a home battery under a tariff
that combines an energy charge with a tiered peak power charge based on the
average of the $N$ largest daily peaks. The tiered charge makes the problem a
mixed-integer linear program, which we solve either directly or by enumerating
tier assignments and solving a sequence of linear programs; both are fast and
reliable. On a year of real data from a home in Trondheim, MPC with simple
forecasts comes within $1.7\%$ of the prescient bound and saves nearly three
times as much as the best rule-based policy we consider.

The formulation is not specific to this tariff or to chemical batteries. The
same MILP applies to any piecewise-constant peak charge and to other storage
assets, including thermal storage and electric-vehicle batteries under
vehicle-to-grid operation.

\section*{Acknowledgments}

David Pérez-Piñeiro was supported by the Research Council of Norway and
HighEFF, an eight-year Research Centre operating under the FME scheme
(Centre for Environment-friendly Energy Research, Grant No. 257632). Stephen Boyd was
partially supported by ACCESS (AI Chip Center for Emerging Smart Systems),
sponsored by InnoHK funding, Hong Kong SAR, and by Office of Naval Research
grant N00014-22-1-2121.

\clearpage
\appendix

\section{Sensitivity analysis}
\label{s-sensitivity}
We study the sensitivity of MPC performance to three design choices: the
planning horizon $H$, the peak power parameter $N$, and the forecast method.
All other parameters are held at their default values from
\S\ref{s-MPC}.

\paragraph{MPC horizon.}
Table~\ref{t-horizon} shows MPC cost as a function of the planning horizon $H$.
A horizon of one day ($H=24$) performs poorly because the controller cannot
anticipate peak power charges beyond the current day; the peak power cost
is $3{,}633$ NOK at $H=24$, versus $2{,}289$ at $H = 720$. Performance improves with longer horizons up to $H = 720$ (30 days), which
spans a full billing month. Doubling the horizon to $H = 1440$ does not
improve performance further; the natural planning horizon is one billing
month, since peak power charges reset monthly.

\begin{table}[h]
    \centering
    \caption{MPC cost versus planning horizon $H$. Costs in NOK.}
    \label{t-horizon}
    \begin{tabular}{lcccc}
        \toprule
        $H$ (hours) & Energy & Peak power & Total & Savings \\
        \midrule
        $24$ (1 day)     & $19{,}213$ & $3{,}633$ & $22{,}846$ & $8.8\%$ \\
        $168$ (1 week)   & $19{,}234$ & $2{,}618$ & $21{,}852$ & $12.8\%$ \\
        $360$ (2 weeks)  & $19{,}280$ & $2{,}513$ & $21{,}793$ & $13.0\%$ \\
        $720$ (30 days)  & $19{,}279$ & $2{,}289$ & $21{,}568$ & $13.9\%$ \\
        $1440$ (60 days) & $19{,}352$ & $2{,}289$ & $21{,}641$ & $13.6\%$ \\
        \bottomrule
    \end{tabular}
\end{table}

\paragraph{Peak power parameter $N$.}
The tariff defines the monthly peak average $z_k$ as the mean of the $N$
largest daily maxima. Table~\ref{t-peak-N} compares $N=1$ and $N=3$.
The difference is negligible (5 NOK over the year). With $N=1$, the
controller optimizes against the single worst day in each month; with
$N=3$ (the actual tariff definition), it can tolerate occasional peaks
as long as their average stays low. The cost is nearly the same either way,
because a policy that limits the largest daily peak also limits the average
of the largest three.

\begin{table}[h]
    \centering
    \caption{MPC cost versus peak power parameter $N$. Costs in NOK.}
    \label{t-peak-N}
    \begin{tabular}{lcccc}
        \toprule
        $N$ & Energy & Peak power & Total & Savings \\
        \midrule
        $1$ & $19{,}274$ & $2{,}289$ & $21{,}563$ & $13.9\%$ \\
        $3$ & $19{,}279$ & $2{,}289$ & $21{,}568$ & $13.9\%$ \\
        \bottomrule
    \end{tabular}
\end{table}

\paragraph{Forecast method.}
Table~\ref{t-forecast} compares MPC using the full baseline+AR forecast to a
baseline-only variant that omits the autoregressive residual component. The
improvement from the AR correction is modest (90 NOK, or 0.4\%) and comes
entirely from reduced peak power charges ($2{,}289$ versus $2{,}394$). The
AR model helps the controller anticipate near-term load peaks, which
improves peak shaving; the seasonal baseline alone captures the price
patterns that drive energy arbitrage.

\begin{table}[h]
    \centering
    \caption{MPC cost versus forecast method. Costs in NOK.}
    \label{t-forecast}
    \begin{tabular}{lcccc}
        \toprule
        Forecast & Energy & Peak power & Total & Savings \\
        \midrule
        Baseline only  & $19{,}264$ & $2{,}394$ & $21{,}658$ & $13.5\%$ \\
        Baseline + AR  & $19{,}279$ & $2{,}289$ & $21{,}568$ & $13.9\%$ \\
        \bottomrule
    \end{tabular}
\end{table}

\clearpage
\printbibliography

\end{document}